\documentclass[12pt]{article}

\usepackage{amsmath}
\usepackage{theorem}
\usepackage{amssymb}
\usepackage{latexsym}
\usepackage{a4wide}
\usepackage{longtable}
\usepackage{epic}

\newtheorem{thm}{Theorem}

\newtheorem{prop}[thm]{Proposition}

{\theorembodyfont{\rmfamily} }
{\theorembodyfont{\rmfamily} \newtheorem{exa}[thm]{Example}}

\newcommand{\C}{\mathbb{C}}

\newcommand{\Z}{\mathbb{Z}}

\newcommand{\ad}{\mathrm{\mathop{ad}}}

\newcommand{\mf}[1]{\mathfrak{#1}}
\newcommand{\ssl}{\mathfrak{\mathop{sl}}}
\newcommand{\g}{\mathfrak{g}}
\newcommand{\myl}{\mathfrak{l}}
\newcommand{\hh}{\mathfrak{h}}

\newcommand{\lab}[1]{\footnotesize{#1}}

\begin{document}

\title{Computations with nilpotent orbits in {\sf SLA}}
\author{Willem A. de Graaf\\
Dipartimento di Matematica\\
Universit\`{a} di Trento\\
Italy}
\date{}
\maketitle

\begin{abstract}
We report on some computations with nilpotent orbits in simple Lie
algebras of exceptional type within the {\sf SLA} package of {\sf GAP}4.
Concerning reachable nilpotent orbits our computations firstly confirm the 
classification of such orbits in Lie algebras of exceptional type 
by Elashvili and Gr\'elaud, secondly they answer a question by Panyushev, 
and thirdly they show in what way a recent result of Yakimova for the
Lie algebras of classical type extends to the exceptional types.
The second topic of this note concerns abelianizations of centralizers of
nilpotent elements. We give tables with their dimensions.
\end{abstract}

\section{Introduction}

In the theory of the nilpotent orbits of simple Lie algebras there appears to be 
a dichotomy between the Lie algebras of classical type and those of exceptional type:
statements for the nilpotent orbits in the Lie algebras of classical type usually
have an elegant proof, whereas the same statements for the exceptional types require
detailed and ``dirty'' calculations (with or without making use of a computer).
This can already be seen in the classification of the nilpotent orbits itself:
for the classical types they are classified in terms of partitions that directly
correspond to the Jordan blocks of a representative, while for the exceptional types
completely different methods have to be used. Another example is the proof of Elashvili's
conjecture: this was proved for the classical types by Yakimova (\cite{yakimova}) and
checked by computer for the exceptional types in \cite{gra14}. Leter Charbonnel and Moreau
found a more uniform proof (\cite{chamo}), which, however, for some cases in the exceptional
types still required computer calculations. \par
This paper is devoted to some more instances of this pattern. In the first half we 
present some results, obtained by computer calculations, on {\em reachable} nilpotent orbits
of simple Lie algebras of exceptional type. They extend results obtained for the classical
types (\cite{panyushev4}, \cite{yakimova2}). In the second half we present computational
data on the algebras $\g_e/[\g_e,\g_e]$. This data has been used by Premet and Topley 
(\cite{pretop}) to extend some of their results, related to finite $W$-algebras, for the
classical types to the exceptional types. (Here the use of ``extend'' includes the possibility
that the statements are slightly changed, or that exceptions are introduced.) \par
All results given here have been obtained by computation - in particular no proofs are
given (the algorithms used have been described elsewhere, cf. \cite{gra14}). This paper is 
intended to serve two purposes: firstly to publish the results of our computations for
future reference, and secondly to advertise the package {\sf SLA} with which the computations
have been done. 

\subsection{Preliminaries on nilpotent orbits}

Now we introduce some notation and recall some definitions and facts on nilpotent orbits.
For more background information we refer to the the book by Collingwood and McGovern 
(\cite{colmcgov}).

Let $\g$ be a simple Lie algebra over $\C$ (or over an algebraically
closed field of characteristic 0). Let $G$ denote the adjoint group of $\g$.
By the Jacobson-Morozov theorem a nilpotent $e\in \g$ lies in an
$\ssl_2$-triple $(h,e,f)$ (where $[e,f]=h$, $[h,e]=2e$, $[h,f]=-2f$).
Let $\hh\subset\g$ be a Cartan subalgebra containing $h$. Let $\Phi$ be the root
system of $\g$ with respect to $\hh$, with basis of simple roots $\Delta = \{\alpha_1,\ldots,
\alpha_\ell\}$. Let $W$ denote the Weyl group of $\Phi$. Then, possibly after replacing $h$ 
with a $W$-conjugate, we have $\alpha_i(h) \in \{0,1,2\}$. The Dynkin diagram of $\Delta$,
where the node corresponding to $\alpha_i$ is labeled $\alpha_i(h)$, is called a
weighted Dynkin diagram. It uniquely determines the orbit $Ge$.

For $e\in \g$ we denote its centraliser
in $\g$ by $\g_e$. In \cite{panyushev4} an $e$ in $\g$ is defined to be
{\em reachable} if $e \in [\g_e,\g_e]$. Such an element has to be nilpotent.
In \cite{elashgre}, Elashvili and Gr\'elaud gave a classification of reachable 
elements in $\g$ (in that paper such elements are called {\em compact}, 
in analogy with \cite{blabry}).  

By the adjoint representation the subalgebra spanned by an $\ssl_2$-triple $(h,e,f)$ acts
on $\g$. Since the eigenvalues of $\ad h$ are integers, we get a grading
$$ \g = \bigoplus_{k\in\Z} \g(k)$$
where $\g(k) = \{ x \in \g\mid [h,x]=kx\}$. Now set $\g(k)_e = \g(k)\cap 
\g_e$, and let $\g(\geq 1)_e$ denote the subalgebra spanned by all
$\g(k)_e$, $k\geq 1$. 

Let $\mf{p}\subset \g$ be a parabolic subalgebra, with Levi decomposition $\mf{p} = \myl\oplus
\mf{n}$, where $\mf{n}$ is the nilradical, and $\myl$ is reductive. Let $L\subset G$ be the
connected subgroup of $G$ with Lie algebra $\myl$. Let $Le'$ be a nilpotent orbit in $\myl$.
Then Lusztig and Spaltenstein (\cite{lusp}) have shown that there is a unique nilpotent
orbit $Ge\subset \g$ such that $Ge\cap (Le'\oplus \mf{n})$ is dense in the latter. 
The orbit $Ge$ is said to be {\em induced} from the orbit $Le'$. Nilpotent orbits which are
not induced are called {\em rigid}.

Let $n$ be a non-negative integer. The irreducible components of the locally closed set 
$$A^n = \{ x \in \g \mid \dim Gx = n\}$$
are called {\em sheets} of $\g$ (see \cite{borho}, \cite{borhokraft}). A sheet is a $G$-stable
subset containing a {\em unique} nilpotent orbit. The converse is not true: sheets in general
are not disjoint, and different sheets may contain the same nilpotent orbit. The sheets of 
$\g$ are in bijection with ($G$-classes of) pairs $(\myl, Le')$, where $\myl$ is a 
Levi subalgebra, and $Le'$ is a rigid nilpotent orbit in $\myl$. The {\em rank} of the sheet
correponding to the pair $(\myl, Le')$ is defined to be the dimension of the centre
of $\myl$.

\subsection{Main results}

Panyushev (\cite{panyushev4}) showed that, for $\g$ of type $A_n$,
$e$ is reachable if and only if $\g(\geq 1)_e$ is generated as Lie algebra
by $\g(1)_e$. Here we call this the {\em Panyushev property} of $\g$.
In \cite{panyushev4} it is stated that this property also holds for the 
other classical types and the question is posed whether it holds for the
exceptional types. In \cite{yakimova2} a proof is given that the Panyushev
property holds in types $B_n$, $C_n$, $D_n$. 
Our computations confirm that the Panyushev property holds also for
the Lie algebras of exceptional type.

Yakimova (\cite{yakimova2}) studied the stronger condition $\g_e = 
[\g_e,\g_e]$. For the purposes of this paper we call elements $e$
satisfying this condition {\em strongly reachable}. She showed that
for $\g$ of classical type, $e$ is strongly reachable if and only if
the nilpotent orbit of $e$ is rigid. Furthermore, 
this is shown to fail for $\g$ of exceptional type. As a result of our
calculations (Section \ref{sec:reach})
we find all rigid nilpotent orbits whose representatives are
not strongly reachable. From this we conclude that $e$ is strongly 
reachable if and only if $e$ is both reachable and rigid.  We note that 
one direction of this statement can be shown in a uniform way for all 
$\g$: if $e$ is strongly reachable then it is reachable, but also rigid by
\cite{yakimova2}, Proposition 11. The converse for exceptional types follows
from our calculations in two ways. Firstly we compute the list of all
strongly reachable orbits and the list of all nilpotent orbits that are
reachable and rigid, and find that they are the same. Second, the Panyushev 
property, which we checked by computation for the exceptional types, 
also implies the statement. (For the classical types we have, of course,
the stronger theorem from \cite{yakimova2}.)

Let $e\in \g$ be a nilpotent element lying in an $\ssl_2$-triple $(h,e,f)$.
In Section \ref{sec:Qe} we consider the quotient $\mf{c}_e = \g_e/[\g_e,\g_e]$, on which
$\ad h$ acts with non-negative eigenvalues. For each $\g$ of exceptional type, we give a table 
listing the dimension of $\mf{c}_e$, where $e$ runs over a set of representatives of the
nilpotent orbits in $\g$, as well as the 
eigenvalues of $\ad h$, acting on $\mf{c}_e$, with their multiplicities. Among other things, 
these tables show that if $Ge$ is an induced nilpotent orbit that lies in a unique
sheet, then the rank of that sheet almost always equals $\dim \mf{c}_e$. The six exceptions
are explicitly listed (Proposition \ref{prop1}). 
Premet and Topley (\cite{pretop}) proved that this holds without 
exceptions for the classical Lie algebras. Proposition \ref{prop1}, as well as the tables of
Section \ref{sec:Qe}, are used in \cite{pretop} for showing that $U(\g,e)^{\mathrm{ab}}$
(the abelianization of a finite $W$-algebra $U(\g,e)$) is isomorphic to a polynomial ring
(with the same six exceptions as Proposition \ref{prop1}).

\subsection{The {\sf SLA} package}

The {\sf SLA} package (\cite{sla} - the acronym stands for Simple Lie Algebras), is written in 
the language of the computer algebra system {\sf GAP}4 (\cite{gap4}), and it can be freely 
downloaded. As the name indicates it has functionality for working with (semi-) 
simple Lie algebras.
Currently there are three main areas which are touched upon by the package: nilpotent
orbits in semisimple Lie algebras, nilpotent orbits of $\theta$-groups, and semisimple subalgebras
of semisimple Lie algebras. 

For the computations underpinning the results of this paper we have used the functionality
for the nilpotent orbits in simple Lie algebras. In particular the package
contains the classification of such orbits. Using this it is 
straightforward to approach the above questions by computational means.
Indeed, for a nilpotent orbit the system easily computes a representative
$e$, and a corresponding $\ssl_2$-triple $(h,e,f)$.
Then using functions present in {\sf GAP}4 we can compute the centralizer,
$\g_e$, and its derived subalgebra, and check whether $e$ lies in it.
This gives us the list of reachable nilpotent orbits. Secondly, a similar
procedure yields the list of strongly reachable orbits.
Thirdly, {\sf SLA} has
a function for computing the grading corresponding to an $\ssl_2$-triple.
With that it is straightforward to check whether $\g(\geq 1)_e$ is generated 
by $\g(1)_e$. Finally, a small procedure for constructing the space $\mf{c}_e = 
\g_e /[\g_e,\g_e]$ and the action of $h$ on it, is easily written. Using that
the results of Section \ref{sec:Qe} are obtained.

{\bf Acknowledgement:} I thank Alexander Elashvili for suggesting the
topic of reachable nilpotent orbits. I thank Alexander Premet for 
suggesting the computations relative to Section \ref{sec:Qe}.

\section{Reachable nilpotent elements in the Lie algebras of exceptional
type}\label{sec:reach}

Tables \ref{tab:rigidE6}, \ref{tab:rigidE7}, \ref{tab:rigidE8}, 
\ref{tab:rigidF4}, and \ref{tab:rigidG2} contain the nilpotent orbits
that by our calculations are reachable. The content of the tables is as
follows. The first column has the label of the orbit, and the second
column the weighted Dynkin diagram. The third and fourth columns contain
a $\times$ if the orbit is, respectively, strongly reachable and rigid.
We note that the classification of rigid nilpotent orbits is known
(see \cite{elashvili}, \cite{elasgra}). 

\setlongtables

\begin{longtable}{|l|c|c|c|}
\caption{Reachable nilpotent orbits in $E_6$.}\label{tab:rigidE6}
\endfirsthead
\hline
\multicolumn{4}{|l|}{\small\slshape Reachable nilpotent orbits in $E_6$.} \\
\hline
\endhead
\hline
\endfoot
\endlastfoot

\hline

label & weighted Dynkin diagram & Strong & Rigid\\
\hline

$A_1$ & $0~~~~0~~~~\overset{\text{\normalsize 1}}{0}~~~~0~~~~0$ 
& $\times$ & $\times$ \\

$2A_1$ & $1~~~~0~~~~\overset{\text{\normalsize 0}}{0}~~~~0~~~~1$ 
&&\\

$3A_1$ & $0~~~~0~~~~\overset{\text{\normalsize 0}}{1}~~~~0~~~~0$ 
& $\times$ & $\times$ \\

$A_2+A_1$ & $1~~~~0~~~~\overset{\text{\normalsize 1}}{0}~~~~0~~~~1$ 
&&\\

$A_2+2A1$ & $0~~~~1~~~~\overset{\text{\normalsize 0}}{0}~~~~1~~~~0$ 
&&\\

$2A_2+A_1$ & $1~~~~0~~~~\overset{\text{\normalsize 0}}{1}~~~~0~~~~1$ 
& $\times$ & $\times$ \\

\hline

\end{longtable}

\begin{longtable}{|l|c|c|c|}
\caption{Reachable nilpotent orbits in $E_7$.}\label{tab:rigidE7}
\endfirsthead
\hline
\multicolumn{4}{|l|}{\small\slshape Reachable nilpotent orbits in $E_7$.} \\
\hline
\endhead
\hline
\endfoot
\endlastfoot

\hline
label & weighted Dynkin diagram & Strong & Rigid\\
\hline

$A_1$ & $1~~~~0~~~~\overset{\text{\normalsize 0}}{0}~~~~0~~~~0~~~~0$ 
& $\times$ & $\times$\\

$2A_1$ & $0~~~~0~~~~\overset{\text{\normalsize 0}}{0}~~~~0~~~~1~~~~0$ 
& $\times$ & $\times$\\

$(3A_1)'$ & $0~~~~1~~~~\overset{\text{\normalsize 0}}{0}~~~~0~~~~0~~~~0$ 
& $\times$ & $\times$\\

$4A_1$ & $0~~~~0~~~~\overset{\text{\normalsize 1}}{0}~~~~0~~~~0~~~~1$ 
& $\times$ & $\times$\\

$A_2+A_1$ & $1~~~~0~~~~\overset{\text{\normalsize 0}}{0}~~~~0~~~~1~~~~0$ 
&&\\

$A_2+2A_1$ & $0~~~~0~~~~\overset{\text{\normalsize 0}}{1}~~~~0~~~~0~~~~0$ 
& $\times$ & $\times$\\

$2A_2+A_1$ & $0~~~~1~~~~\overset{\text{\normalsize 0}}{0}~~~~0~~~~1~~~~0$ 
& $\times$ & $\times$\\

$A_4+A_1$ & $1~~~~0~~~~\overset{\text{\normalsize 0}}{1}~~~~0~~~~1~~~~0$ 
&&\\

\hline

\end{longtable}

\begin{longtable}{|l|c|c|c|}
\caption{Reachable nilpotent orbits in $E_8$.}\label{tab:rigidE8}
\endfirsthead
\hline
\multicolumn{4}{|c|}{\small\slshape Reachable nilpotent orbits in $E_8$.} \\
\hline
\endhead
\hline
\endfoot
\endlastfoot

\hline
label & weighted Dynkin diagram & Strong & Rigid\\
\hline

$A_1$ &
$0~~~~0~~~~\overset{\text{\normalsize 0}}{0}~~~~0~~~~0~~~~0~~~~1$ 
& $\times$ & $\times$\\

$2A_1$ &
$1~~~~0~~~~\overset{\text{\normalsize 0}}{0}~~~~0~~~~0~~~~0~~~~0$
& $\times$ & $\times$\\

$3A_1$ &
$0~~~~0~~~~\overset{\text{\normalsize 0}}{0}~~~~0~~~~0~~~~1~~~~0$
& $\times$ & $\times$\\

$4A_1$ &
$0~~~~0~~~~\overset{\text{\normalsize 1}}{0}~~~~0~~~~0~~~~0~~~~0$
& $\times$ & $\times$\\

$A_2+A_1$ &
$1~~~~0~~~~\overset{\text{\normalsize 0}}{0}~~~~0~~~~0~~~~0~~~~1$
& $\times$ & $\times$\\

$A_2+2A_1$ &
$0~~~~0~~~~\overset{\text{\normalsize 0}}{0}~~~~0~~~~1~~~~0~~~~0$
& $\times$ & $\times$\\

$A_2+3A_1$ &
$0~~~~1~~~~\overset{\text{\normalsize 0}}{0}~~~~0~~~~0~~~~0~~~~0$
& $\times$ & $\times$\\

$2A_2+A_1$ &
$1~~~~0~~~~\overset{\text{\normalsize 0}}{0}~~~~0~~~~0~~~~1~~~~0$
& $\times$ & $\times$\\

$A_4+A_1$ &
$1~~~~0~~~~\overset{\text{\normalsize 0}}{0}~~~~0~~~~1~~~~0~~~~1$
&&\\

$2A_2+2A_1$ &
$0~~~~0~~~~\overset{\text{\normalsize 0}}{0}~~~~1~~~~0~~~~0~~~~0$
& $\times$ & $\times$\\

$(A_3+2A_1)''$ &
$0~~~~1~~~~\overset{\text{\normalsize 0}}{0}~~~~0~~~~0~~~~0~~~~1$
& $\times$ & $\times$\\

$D_4(a_1)+A_1$ &
$0~~~~0~~~~\overset{\text{\normalsize 1}}{0}~~~~0~~~~0~~~~1~~~~0$
& $\times$ & $\times$\\

$A_3+A_2+A_1$ &
$0~~~~0~~~~\overset{\text{\normalsize 0}}{1}~~~~0~~~~0~~~~0~~~~0$
& $\times$ & $\times$\\

$2A_3$ &
$1~~~~0~~~~\overset{\text{\normalsize 0}}{0}~~~~1~~~~0~~~~0~~~~0$ 
& $\times$ & $\times$\\

$A_4+2A_1$ &
$0~~~~0~~~~\overset{\text{\normalsize 0}}{1}~~~~0~~~~0~~~~0~~~~1$ 
&&\\

$A_4+A_3$ &
$0~~~~0~~~~\overset{\text{\normalsize 0}}{1}~~~~0~~~~0~~~~1~~~~0$ 
& $\times$ & $\times$\\

\hline

\end{longtable}

\begin{longtable}{|l|c|c|c|}
\caption{Reachable nilpotent orbits in $F_4$.}\label{tab:rigidF4}
\endfirsthead
\hline
\multicolumn{4}{|l|}{\small\slshape Reachable nilpotent orbits in $F_4$.} \\
\hline
\endhead
\hline
\endfoot
\endlastfoot

\hline

label & weighted Dynkin diagram & Strong & Rigid\\
\hline
 &
\begin{picture}(20,7)
  \put(-20,0){\circle{6}}
  \put(0,0){\circle{6}}
  \put(20,0){\circle{6}}
  \put(40,0){\circle{6}}
  \put(-17,0){\line(1,0){14}}
  \put(2,2){\line(1,0){16}}
  \put(2,-2){\line(1,0){16}}
\put(5,-3){$>$}
  \put(23,0){\line(1,0){14}}
\end{picture} & &

\\
\hline

$A_1$ & 1~~~0~~~0~~~0
& $\times$ & $\times$\\

$\widetilde{A}_1$ & 0~~~0~~~0~~~1
& $\times$ & $\times$\\

$A_1+\widetilde{A}_1$ & 0~~~1~~~0~~~0
& $\times$ & $\times$\\

$A_2+\widetilde{A}_1$ & 0~~~0~~~1~~~0
& $\times$ & $\times$\\

\hline

\end{longtable}

\begin{longtable}{|l|c|c|c|}
\caption{Reachable nilpotent orbits in $G_2$.}\label{tab:rigidG2}
\endfirsthead
\hline
\multicolumn{4}{|l|}{\small\slshape Reachable nilpotent orbits in $G_2$.} \\
\hline
\endhead
\hline
\endfoot
\endlastfoot

\hline

label & weighted Dynkin diagram & Strong & Rigid\\

&

  \begin{picture}(10,7)
  \put(-10,0){\circle{6}}
  \put(22,0){\circle{6}}
\put(3,-3){$>$}
  \put(-8,-2){\line(1,0){28}}
  \put(-7,0){\line(1,0){26}}
  \put(-8,2){\line(1,0){28}}
\end{picture} &&                   \\
\hline

\hline

$\widetilde{A}_1$ & 1~~0 
& $\times$ & $\times$
\\

\hline

\end{longtable}

We make the following comments.

\begin{itemize}
\item Here the reachable elements are exactly the same as in the
paper of Elashvili and Gr\'elaud. Therefore our calculations confirm their
result.
\item The rigid nilpotent orbits that are not strongly reachable are
\begin{itemize}
\item in type $E_7$: $(A_3+A_1)'$ $(41,40)$,
\item in type $E_8$: $A_3+A_1$ $(84,83)$, $D_5(a_1)+A_2$ $(46,45)$, 
$A_5+A_1$ $(46,45)$,
\item in type $F_4$: $\widetilde{A}_2+A_1$ $(16,15)$,
\item in type $G_2$: $A_1$ $(6,5)$.
\end{itemize}
Here the pair of integers in brackets is $(\dim \g_e,\dim [\g_e,\g_e])$.

\item In type $E_6$ all rigid orbits are strongly reachable. Hence 
in this type the situation is the same as for the classical types: 
$e$ is strongly reachable if
and only if the orbit of $e$ is rigid.
\item The last two columns of all tables are equal. This shows that
for the exceptional types the following theorem holds: $e$ is strongly
reachable if and only if $e$ is both reachable and rigid. 
\item This last statement also follows from the Panyushev property. 
Indeed, $e$ rigid implies that $\g(0)_e$ is semisimple, so
$[\g(0)_e,\g(0)_e]=\g(0)_e$. Furthermore, $[\g(0)_e,\g(1)_e]= 
\g(1)_e$ by \cite{yakimova2}, Lemma 8 (where this is shown to hold for 
all nilpotent $e$).
By the Panyushev property this implies that $[\g_e,\g_e]= \g_e$.
\item We see that for all nilpotent orbits that are rigid but not
stronly reachable the codimension of $[\g_e,\g_e]$ in $\g_e$ is 1.
Since a rigid orbit is reachable if and only if it is strongly reachable,
we get that in all those cases $e$ spans the quotient $\g_e/[\g_e,\g_e]$.
\end{itemize}

\begin{exa}
Let us consider the nilpotent orbit in the Lie algebra of type $E_7$ with
label $A_3+A_2$. This orbit is not reachable. It has a representative
with diagram

\begin{picture}(5,35) 
\put(3,10){\circle{6}} 
\put(5,15){\lab{29}} 
\put(23,10){\circle{6}}
\put(25,15){\lab{32}}
\put(6,10){\line(1,0){14}}
\put(43,10){\circle{6}}
\put(26,10){\line(1,0){14}}
\put(45,15){\lab{31}}
\put(63,10){\circle{6}}
\put(65,15){\lab{27}}
\put(83,10){\circle{6}}
\put(85,15){\lab{30}}
\put(66,10){\line(1,0){14}}
\end{picture}

This means that the representative is $e=x_{29}+x_{32}+x_{31}+x_{27}+x_{30}$,
where $x_i$ denotes the root vector corresponding to the $i$-th positive 
root (enumeration as in {\sf GAP}4, cf. \cite{gra14}). Furthermore,
the Dynkin diagram of these roots is as shown above. This representative
is stored in the package {\sf SLA}.

Now, if the orbit were reachable then $e \in [\g_e,\g_e]\cap \g(2)$.
Using the {\sf SLA} package we can easiliy compute the latter space:

\begin{verbatim}
gap> L:= SimpleLieAlgebra("E",7,Rationals);;
gap> o:= NilpotentOrbits(L);;
gap> sl2:=SL2Triple( o[19] );
[ (2)*v.90+(3)*v.92+(2)*v.93+(3)*v.94+(4)*v.95, (6)*v.127+(9)*v.128+(12)*v.129
+(18)*v.130+(14)*v.131+(10)*v.132+(5)*v.133, v.27+v.29+v.30+v.31+v.32 ]
gap> g:= SL2Grading( L, sl2[2] );;
gap> g2:= Subspace( L, g[1][2] );;
gap> der:= LieDerivedSubalgebra(LieCentralizer(L,Subalgebra(L,[sl2[3]])));
<Lie algebra of dimension 33 over Rationals>
gap> BasisVectors( Basis( Intersection( g2, der ) ) );
[ v.18, v.23+(-1)*v.24+v.28, v.24+(-1)*v.25+(-1)*v.28, 
v.27+(-1)*v.29+v.30+(-1)*v.31+(-1)*v.32, v.33+(-1)*v.36+v.37, 
v.34+(-1)*v.36+v.37, v.39 ]
\end{verbatim}

First we make some comments on the above computation. The $\ssl_2$-triple
comes ordered as $(f,h,e)$. So the second element is the neutral element,
and the third element is the nil-positive element, i.e., the representative,
which is as indicated above. So the second element defines the grading,
which we compute with {\tt SL2Grading}. In the subsequent line the
subspace $\g(2)$ is defined, followed by $[\g_e,\g_e]$. Finally a basis
of the intersection is computed.

We see that one of the basis vectors of the intersection is 

\begin{verbatim}
v.27+(-1)*v.29+v.30+(-1)*v.31+(-1)*v.32.
\end{verbatim}

So we see that $e=(x_{29}+x_{32}+x_{31})+(x_{27}+x_{30})$ does not lie
in $[\g_e,\g_e]$ but $(x_{29}+x_{32}+x_{31})-(x_{27}+x_{30})$ does!
\end{exa}

\section{The quotients $\mf{c}_e$}\label{sec:Qe}

Let $\g$ be a simple Lie algebra, and $e$ a representative of a nilpotent orbit
lying in an $\ssl_2$-triple $(f,h,e)$. Note that the centraliser $\g_e$ is spanned
by $\ad h$-eigenvectors with non-negative eigenvalues. We consider the quotient
$\mf{c}_e = \g_e/[\g_e,\g_e]$, on which $\ad h$ also acts. In Tables \ref{tab:E6},
 \ref{tab:E7},  \ref{tab:E8},  \ref{tab:F4},  \ref{tab:G2} we have listed 
the dimension of $\mf{c}_e$ for all nilpotent orbits of the Lie algebras of
exceptional type, along with the eigenvalues of $\ad h$ acting on $\mf{c}_e$.
Comparing these tables and the tables of \cite{elasgra} we get the following two results.

\begin{prop}\label{prop1}
Let $\g$ be a simple Lie algebra of exceptional type. Let $e\in \g$ be a representative
of an induced nilpotent orbit lying in a unique sheet. Then the rank of that sheet is equal to
the dimension of $\g_e/[\g_e,\g_e]$, except the cases listed in Table \ref{tab:prop1}.
\end{prop}

\begin{longtable}{|l|l|l|c|c|}
\caption{Table of exceptions to Proposition \ref{prop1}.}\label{tab:prop1}
\endfirsthead
\hline
\multicolumn{5}{|l|}{Exceptions.} \\
\hline
\endhead
\hline
\endfoot
\endlastfoot

\hline

$\g$ & label & weighted Dynkin diagram & rank & $\dim \g_e/[\g_e,\g_e]$ \\
\hline

$E_6$ & $A_3+A_1$ & $0~~~~1~~~~\overset{\text{\normalsize 1}}{0}~~~~1~~~~0$ 
& 1 & 2 \\
$E_7$ & $D_6(a_2)$ & $0~~~~1~~~~\overset{\text{\normalsize 1}}{0}~~~~1~~~~0~~~~2$
& 2 & 3\\
$E_8$ & $D_6(a_2)$ & $0~~~~1~~~~\overset{\text{\normalsize 1}}{0}~~~~0~~~~0~~~~1~~~~0$
& 1 & 3\\
$E_8$ & $E_6(a_3)+A_1$ & $1~~~~0~~~~\overset{\text{\normalsize 0}}{0}~~~~1~~~~0~~~~1~~~~0$
& 1 & 3\\
$E_8$ & $E_7(a_2)$ & $0~~~~1~~~~\overset{\text{\normalsize 1}}{0}~~~~1~~~~0~~~~2~~~~2$
& 3 & 4\\
$F_4$ & $C_3(a_1)$ &  1~~~0~~~1~~~0 & 1 & 3\\
\hline
\end{longtable}

\begin{prop}\label{prop2}
Let $\g$ be a simple Lie algebra of exceptional type, and let $e\in \g$ be a nilpotent
orbit that lies in more than one sheet. Then the maximal rank of such a sheet
is strictly smaller than $\dim \g_e/[\g_e,\g_e]$.
\end{prop}

\begin{longtable}{|l|c|c|l|}
\caption{Nilpotent orbits in the Lie algebra of type $E_6$.}\label{tab:E6}
\endfirsthead
\hline
\multicolumn{4}{|l|}{\small\slshape Nilpotent orbits in type $E_6$} \\ 
\hline
\endhead
\hline
\endfoot
\endlastfoot

\hline
label & weighted Dynkin diagram & $\dim \mf{c}_e$ & $h$-weights \\

&  

\begin{picture}(40,25)
  \put(-20,0){\circle{6}}
  \put(0,0){\circle{6}}
  \put(20,0){\circle{6}}
  \put(40,0){\circle{6}}
  \put(60,0){\circle{6}}
  \put(20,20){\circle{6}}
  \put(-17,0){\line(1,0){14}}
  \put(3,0){\line(1,0){14}}
  \put(23,0){\line(1,0){14}}
  \put(43,0){\line(1,0){14}}
  \put(20,3){\line(0,1){14}}
\end{picture}                   

& & \\
\hline

$A_1$ & $0~~~~0~~~~\overset{\text{\normalsize 1}}{0}~~~~0~~~~0$ & 
0 & \\

$2A_1$ & $1~~~~0~~~~\overset{\text{\normalsize 0}}{0}~~~~0~~~~1$ & 
1 & 0 \\

$3A_1$ & $0~~~~0~~~~\overset{\text{\normalsize 0}}{1}~~~~0~~~~0$ & 
0 & \\

$A_2$ & $0~~~~0~~~~\overset{\text{\normalsize 2}}{0}~~~~0~~~~0$ & 
1 & 2\\

$A_2+A_1$ & $1~~~~0~~~~\overset{\text{\normalsize 1}}{0}~~~~0~~~~1$ & 
1 & 0 \\

$2A_2$ & $2~~~~0~~~~\overset{\text{\normalsize 0}}{0}~~~~0~~~~2$ & 
2 & 2,4\\

$A_2+2A_1$ & $0~~~~1~~~~\overset{\text{\normalsize 0}}{0}~~~~1~~~~0$ & 
1 & 0\\

$A_3$ & $1~~~~0~~~~\overset{\text{\normalsize 2}}{0}~~~~0~~~~1$ & 
2 & 0, 2\\

$2A_2+A_1$ & $1~~~~0~~~~\overset{\text{\normalsize 0}}{1}~~~~0~~~~1$ & 
0 & \\

$A_3+A_1$ & $0~~~~1~~~~\overset{\text{\normalsize 1}}{0}~~~~1~~~~0$ & 
2 & 0,2\\

$D_4(a_1)$ & $0~~~~0~~~~\overset{\text{\normalsize 0}}{2}~~~~0~~~~0$ & 
5 & 0,0,2,2,2\\

$A_4$ & $2~~~~0~~~~\overset{\text{\normalsize 2}}{0}~~~~0~~~~2$ & 
3 & 0,2,6\\

$D_4$ & $0~~~~0~~~~\overset{\text{\normalsize 2}}{2}~~~~0~~~~0$ & 
2 & 2,10\\

$A_4+A_1$ & $1~~~~1~~~~\overset{\text{\normalsize 1}}{0}~~~~1~~~~1$ & 
2 & 0,2 \\

$A_5$ & $2~~~~1~~~~\overset{\text{\normalsize 1}}{0}~~~~1~~~~2$ & 
2 & 2,4\\

$D_5(a_1)$ & $1~~~~1~~~~\overset{\text{\normalsize 2}}{0}~~~~1~~~~1$ & 
3 & 0,2,4\\

$E_6(a_3)$ & $2~~~~0~~~~\overset{\text{\normalsize 0}}{2}~~~~0~~~~2$ & 
5 & 2,2,2,4,4\\

$D_5$ & $2~~~~0~~~~\overset{\text{\normalsize 2}}{2}~~~~0~~~~2$ & 
4 & 0,2,6,10\\

$E_6(a_1)$ & $2~~~~2~~~~\overset{\text{\normalsize 2}}{0}~~~~2~~~~2$ & 
5 & 2,4,6,8,10\\

$E_6$ & $2~~~~2~~~~\overset{\text{\normalsize 2}}{2}~~~~2~~~~2$ & 
6 & 2,8,10,14,16,22\\

\hline
\end{longtable}

\begin{longtable}{|l|c|c|l|}
\caption{Nilpotent orbits in the Lie algebra of type $E_7$.}\label{tab:E7}
\endfirsthead
\hline
\multicolumn{4}{|l|}{\small Nilpotent orbits in $E_7$.} \\ 
\hline
\endhead
\hline
\endfoot
\endlastfoot

\hline
label & weighted Dynkin diagram & $\dim \mf{c}_e$ & $h$-weights \\

&  

\begin{picture}(100,25)
  \put(0,0){\circle{6}}
  \put(20,0){\circle{6}}
  \put(40,0){\circle{6}}
  \put(60,0){\circle{6}}
  \put(80,0){\circle{6}}
  \put(100,0){\circle{6}}
  \put(40,20){\circle{6}}
  \put(3,0){\line(1,0){14}}
  \put(23,0){\line(1,0){14}}
  \put(43,0){\line(1,0){14}}
  \put(63,0){\line(1,0){14}}
  \put(83,0){\line(1,0){14}}
  \put(40,3){\line(0,1){14}}
\end{picture}                   

& & \\
\hline

$A_1$ & $1~~~~0~~~~\overset{\text{\normalsize 0}}{0}~~~~0~~~~0~~~~0$ & 
0 & \\

$2A_1$ & $0~~~~0~~~~\overset{\text{\normalsize 0}}{0}~~~~0~~~~1~~~~0$ & 
0 & \\

$(3A_1)''$ & $0~~~~0~~~~\overset{\text{\normalsize 0}}{0}~~~~0~~~~0~~~~2$ & 
1 & 2 \\

$(3A_1)'$ & $0~~~~1~~~~\overset{\text{\normalsize 0}}{0}~~~~0~~~~0~~~~0$ & 
0 & \\

$A_2$ & $2~~~~0~~~~\overset{\text{\normalsize 0}}{0}~~~~0~~~~0~~~~0$ & 
1 & 2\\

$4A_1$ & $0~~~~0~~~~\overset{\text{\normalsize 1}}{0}~~~~0~~~~0~~~~1$ & 
0 & \\

$A_2+A_1$ & $1~~~~0~~~~\overset{\text{\normalsize 0}}{0}~~~~0~~~~1~~~~0$ & 
1 & 0 \\

$A_2+2A_1$ & $0~~~~0~~~~\overset{\text{\normalsize 0}}{1}~~~~0~~~~0~~~~0$ & 
0 & \\

$A_3$ & $2~~~~0~~~~\overset{\text{\normalsize 0}}{0}~~~~0~~~~1~~~~0$ & 
1 & 2 \\

$2A_2$ & $0~~~~0~~~~\overset{\text{\normalsize 0}}{0}~~~~0~~~~2~~~~0$ & 
1 & 2\\

$A_2+3A_1$ & $0~~~~0~~~~\overset{\text{\normalsize 2}}{0}~~~~0~~~~0~~~~0$ & 
1 & 2 \\

$(A_3+A_1)''$ & $2~~~~0~~~~\overset{\text{\normalsize 0}}{0}~~~~0~~~~0~~~~2$ & 
2 & 2,2\\

$2A_2+A_1$ & $0~~~~1~~~~\overset{\text{\normalsize 0}}{0}~~~~0~~~~1~~~~0$ & 
0 & \\

$(A_3+A_1)'$ & $1~~~~0~~~~\overset{\text{\normalsize 0}}{1}~~~~0~~~~0~~~~0$ & 
1 & 2\\

$D_4(a_1)$ & $0~~~~2~~~~\overset{\text{\normalsize 0}}{0}~~~~0~~~~0~~~~0$ & 
3 & 2,2,2\\

$A_3+2A_1$ & $1~~~~0~~~~\overset{\text{\normalsize 0}}{0}~~~~1~~~~0~~~~1$ & 
1 & 2\\

$D_4$ & $2~~~~2~~~~\overset{\text{\normalsize 0}}{0}~~~~0~~~~0~~~~0$ &
2 & 2,10\\

$D_4(a_1)+A_1$ & $0~~~~1~~~~\overset{\text{\normalsize 1}}{0}~~~~0~~~~0~~~~1$ &
2 & 2,2\\

$A_3+A_2$ & $0~~~~0~~~~\overset{\text{\normalsize 0}}{1}~~~~0~~~~1~~~~0$ & 
2 & 0,2\\

$A_4$ & $2~~~~0~~~~\overset{\text{\normalsize 0}}{0}~~~~0~~~~2~~~~0$ & 
3 & 0,2,6\\

$A_3+A_2+A_1$ & $0~~~~0~~~~\overset{\text{\normalsize 0}}{0}~~~~2~~~~0~~~~0$ & 
1 & 2\\

$(A_5)''$ & $2~~~~0~~~~\overset{\text{\normalsize 0}}{0}~~~~0~~~~2~~~~2$ & 
3 & 2,6,10\\

$D_4+A_1$ & $2~~~~1~~~~\overset{\text{\normalsize 1}}{0}~~~~0~~~~0~~~~1$ &
1 & 2\\

$A_4+A_1$ & $1~~~~0~~~~\overset{\text{\normalsize 0}}{1}~~~~0~~~~1~~~~0$ & 
2 & 0,0\\

$D_5(a_1)$ & $2~~~~0~~~~\overset{\text{\normalsize 0}}{1}~~~~0~~~~1~~~~0$ &
3 & 0,2,4\\

$A_4+A_2$ & $0~~~~0~~~~\overset{\text{\normalsize 0}}{2}~~~~0~~~~0~~~~0$ & 
1 & 2\\

$(A_5)'$ & $1~~~~0~~~~\overset{\text{\normalsize 0}}{1}~~~~0~~~~2~~~~0$ & 
1 & 2\\

$A_5+A_1$ & $1~~~~0~~~~\overset{\text{\normalsize 0}}{1}~~~~0~~~~1~~~~2$ & 
1 & 2\\

$D_5(a_1)+A_1$ & $2~~~~0~~~~\overset{\text{\normalsize 0}}{0}~~~~2~~~~0~~~~0$ & 
2 & 2,2\\

$D_6(a_2)$ & $0~~~~1~~~~\overset{\text{\normalsize 1}}{0}~~~~1~~~~0~~~~2$ & 
3 & 2,2,2\\

$E_6(a_3)$ & $0~~~~2~~~~\overset{\text{\normalsize 0}}{0}~~~~0~~~~2~~~~0$ & 
3 & 2,2,2\\

$D_5$ & $2~~~~2~~~~\overset{\text{\normalsize 0}}{0}~~~~0~~~~2~~~~0$ & 
3 & 2,6,10\\

$E_7(a_5)$ & $0~~~~0~~~~\overset{\text{\normalsize 0}}{2}~~~~0~~~~0~~~~2$ & 
6 & 2,2,2,2,2,2\\

$A_6$ & $0~~~~0~~~~\overset{\text{\normalsize 0}}{2}~~~~0~~~~2~~~~0$ & 
2 & 2,10\\

$D_5+A_1$ & $2~~~~1~~~~\overset{\text{\normalsize 1}}{0}~~~~1~~~~1~~~~0$ & 
3 & 2,2,10\\

$D_6(a_1)$ & $2~~~~1~~~~\overset{\text{\normalsize 1}}{0}~~~~1~~~~0~~~~2$ & 
4 & 2,2,6,10\\

$E_7(a_4)$ & $2~~~~0~~~~\overset{\text{\normalsize 0}}{2}~~~~0~~~~0~~~~2$ & 
4 & 2,2,2,2\\

$D_6$ & $2~~~~1~~~~\overset{\text{\normalsize 1}}{0}~~~~1~~~~2~~~~2$ & 
3 & 2,6,10\\

$E_6(a_1)$ & $2~~~~0~~~~\overset{\text{\normalsize 0}}{2}~~~~0~~~~2~~~~0$ & 
5 & 0,2,4,6,10\\

$E_6$ & $2~~~~2~~~~\overset{\text{\normalsize 0}}{2}~~~~0~~~~2~~~~0$ & 
4 & 2,10,14,22\\

$E_7(a_3)$ & $2~~~~0~~~~\overset{\text{\normalsize 0}}{2}~~~~0~~~~2~~~~2$ & 
6 & 2,2,4,6,8,10\\

$E_7(a_2)$ & $2~~~~2~~~~\overset{\text{\normalsize 2}}{0}~~~~2~~~~0~~~~2$ & 
5 & 2,2,6,10,14\\

$E_7(a_1)$ & $2~~~~2~~~~\overset{\text{\normalsize 2}}{0}~~~~2~~~~2~~~~2$ & 
6 & 2,6,10,10,14,18\\

$E_7$ & $2~~~~2~~~~\overset{\text{\normalsize 2}}{2}~~~~2~~~~2~~~~2$ & 
7 & 2,10,14,18,22,26,34\\

\hline
\end{longtable}

\begin{longtable}{|l|c|c|l|}
\caption{Nilpotent orbits in the Lie algebra of type $E_8$.}\label{tab:E8}
\endfirsthead
\hline
\multicolumn{4}{|l|}{\small Nilpotent orbits in $E_8$.} \\ 
\hline
\endhead
\hline
\endfoot
\endlastfoot

\hline
label & weighted Dynkin diagram & $\dim \mf{c}_e$ & $h$-weights \\

&  

\begin{picture}(120,25)
  \put(0,0){\circle{6}}
  \put(20,0){\circle{6}}
  \put(40,0){\circle{6}}
  \put(60,0){\circle{6}}
  \put(80,0){\circle{6}}
  \put(100,0){\circle{6}}
  \put(120,0){\circle{6}}
  \put(40,20){\circle{6}}
  \put(3,0){\line(1,0){14}}
  \put(23,0){\line(1,0){14}}
  \put(43,0){\line(1,0){14}}
  \put(63,0){\line(1,0){14}}
  \put(83,0){\line(1,0){14}}
  \put(40,3){\line(0,1){14}}
  \put(103,0){\line(1,0){14}}
\end{picture}                   

& & \\
\hline

$A_1$ & $0~~~~0~~~~\overset{\text{\normalsize 0}}{0}~~~~0~~~~0~~~~0~~~~1$ & 
0 &\\

$2A_1$ & $1~~~~0~~~~\overset{\text{\normalsize 0}}{0}~~~~0~~~~0~~~~0~~~~0$ & 
0 &\\

$3A_1$ & $0~~~~0~~~~\overset{\text{\normalsize 0}}{0}~~~~0~~~~0~~~~1~~~~0$ & 
0 &\\

$A_2$ & $0~~~~0~~~~\overset{\text{\normalsize 0}}{0}~~~~0~~~~0~~~~0~~~~2$ & 
1 & 2\\

$4A_1$ & $0~~~~0~~~~\overset{\text{\normalsize 1}}{0}~~~~0~~~~0~~~~0~~~~0$ & 
0 & \\

$A_2+A_1$ & $1~~~~0~~~~\overset{\text{\normalsize 0}}{0}~~~~0~~~~0~~~~0~~~~1$ & 
0 & \\

$A_2+2A_1$ & $0~~~~0~~~~\overset{\text{\normalsize 0}}{0}~~~~0~~~~1~~~~0~~~~0$ & 
0 & \\

$A_3$ & $1~~~~0~~~~\overset{\text{\normalsize 0}}{0}~~~~0~~~~0~~~~0~~~~2$ & 
1 & 2\\

$A_2+3A_1$ & $0~~~~1~~~~\overset{\text{\normalsize 0}}{0}~~~~0~~~~0~~~~0~~~~0$ & 
0 & \\

$2A_2$ & $2~~~~0~~~~\overset{\text{\normalsize 0}}{0}~~~~0~~~~0~~~~0~~~~0$ & 
1 & 2\\

$2A_2+A_1$ & $1~~~~0~~~~\overset{\text{\normalsize 0}}{0}~~~~0~~~~0~~~~1~~~~0$ & 
0 & \\

$A_3+A_1$ & $0~~~~0~~~~\overset{\text{\normalsize 0}}{0}~~~~0~~~~1~~~~0~~~~1$ & 
1 & 2 \\

$D_4(a_1)$ & $0~~~~0~~~~\overset{\text{\normalsize 0}}{0}~~~~0~~~~0~~~~2~~~~0$ & 
3 & 2,2,2\\

$D_4$ & $0~~~~0~~~~\overset{\text{\normalsize 0}}{0}~~~~0~~~~0~~~~2~~~~2$ &
2 & 2,10\\

$2A_2+2A_1$ & $0~~~~0~~~~\overset{\text{\normalsize 0}}{0}~~~~1~~~~0~~~~0~~~~0$ & 
0 & \\

$A_3+2A_1$ & $0~~~~1~~~~\overset{\text{\normalsize 0}}{0}~~~~0~~~~0~~~~0~~~~1$ & 
0 & \\

$D_4(a_1)+A_1$ & $0~~~~0~~~~\overset{\text{\normalsize 1}}{0}~~~~0~~~~0~~~~1~~~~0$ & 
0 & \\

$A_3+A_2$ & $1~~~~0~~~~\overset{\text{\normalsize 0}}{0}~~~~0~~~~1~~~~0~~~~0$ & 
2 & 0,2\\

$A_4$ & $2~~~~0~~~~\overset{\text{\normalsize 0}}{0}~~~~0~~~~0~~~~0~~~~2$ & 
2 & 2,6\\

$A_3+A_2+A_1$ & $0~~~~0~~~~\overset{\text{\normalsize 0}}{1}~~~~0~~~~0~~~~0~~~~0$ & 
0 & \\

$D_4+A_1$ & $0~~~~0~~~~\overset{\text{\normalsize 1}}{0}~~~~0~~~~0~~~~1~~~~2$ & 
1 & 2\\

$D_4(a_1)+A_2$ & $0~~~~0~~~~\overset{\text{\normalsize 2}}{0}~~~~0~~~~0~~~~0~~~~0$ & 
1 & 2\\

$A_4+A_1$ & $1~~~~0~~~~\overset{\text{\normalsize 0}}{0}~~~~0~~~~1~~~~0~~~~1$ & 
1 & 0\\

$2A_3$ & $1~~~~0~~~~\overset{\text{\normalsize 0}}{0}~~~~1~~~~0~~~~0~~~~0$ & 
0 & \\

$D_5(a_1)$ & $1~~~~0~~~~\overset{\text{\normalsize 0}}{0}~~~~0~~~~1~~~~0~~~~2$ & 
2 & 2,4\\

$A_4+2A_1$ & $0~~~~0~~~~\overset{\text{\normalsize 0}}{1}~~~~0~~~~0~~~~0~~~~1$ & 
1 & 0\\

$A_4+A_2$ & $0~~~~0~~~~\overset{\text{\normalsize 0}}{0}~~~~0~~~~2~~~~0~~~~0$ & 
1 & 2 \\

$A_5$ & $2~~~~0~~~~\overset{\text{\normalsize 0}}{0}~~~~0~~~~1~~~~0~~~~1$ & 
1 & 2 \\

$D_5(a_1)+A_1$ & $0~~~~0~~~~\overset{\text{\normalsize 0}}{1}~~~~0~~~~0~~~~0~~~~2$ & 
1 & 2 \\

$A_4+A_2+A_1$ & $0~~~~1~~~~\overset{\text{\normalsize 0}}{0}~~~~0~~~~1~~~~0~~~~0$ & 
1 & 2\\

$D_4+A_2$ & $0~~~~0~~~~\overset{\text{\normalsize 2}}{0}~~~~0~~~~0~~~~0~~~~2$ & 
2 & 2,2\\

$E_6(a_3)$ & $2~~~~0~~~~\overset{\text{\normalsize 0}}{0}~~~~0~~~~0~~~~2~~~~0$ & 
3 & 2,2,2\\

$D_5$ & $2~~~~0~~~~\overset{\text{\normalsize 0}}{0}~~~~0~~~~0~~~~2~~~~2$ & 
3 & 2,6,10\\

$A_4+A_3$ & $0~~~~0~~~~\overset{\text{\normalsize 0}}{1}~~~~0~~~~0~~~~1~~~~0$ & 
0 & \\

$A_5+A_1$ & $1~~~~0~~~~\overset{\text{\normalsize 0}}{1}~~~~0~~~~0~~~~0~~~~1$ & 
1 & 2\\

$D_5(a_1)+A_2$ & $0~~~~1~~~~\overset{\text{\normalsize 0}}{0}~~~~0~~~~1~~~~0~~~~1$ & 
1 & 2\\

$D_6(a_2)$ & $0~~~~1~~~~\overset{\text{\normalsize 1}}{0}~~~~0~~~~0~~~~1~~~~0$ & 
3 & 2,2,2\\

$E_6(a_3)+A_1$ & $1~~~~0~~~~\overset{\text{\normalsize 0}}{0}~~~~1~~~~0~~~~1~~~~0$ & 
3 & 2,2,2\\

$E_7(a_5)$ & $0~~~~0~~~~\overset{\text{\normalsize 0}}{1}~~~~0~~~~1~~~~0~~~~0$ & 
6 & 2,2,2,2,2,2\\

$D_5+A_1$ & $1~~~~0~~~~\overset{\text{\normalsize 0}}{0}~~~~1~~~~0~~~~1~~~~2$ & 
2 & 2,10\\

$E_8(a_7)$ & $0~~~~0~~~~\overset{\text{\normalsize 0}}{0}~~~~2~~~~0~~~~0~~~~0$ & 
10 & 2,2,2,2,2,2,2,2,2,2\\

$A_6$ & $2~~~~0~~~~\overset{\text{\normalsize 0}}{0}~~~~0~~~~2~~~~0~~~~0$ & 
2 & 2,10\\

$D_6(a_1)$ & $0~~~~1~~~~\overset{\text{\normalsize 1}}{0}~~~~0~~~~0~~~~1~~~~2$ & 
3 & 2,2,10\\

$A_6+A_1$ & $1~~~~0~~~~\overset{\text{\normalsize 0}}{1}~~~~0~~~~1~~~~0~~~~0$ & 
1 & 2\\

$E_7(a_4)$ & $0~~~~0~~~~\overset{\text{\normalsize 0}}{1}~~~~0~~~~1~~~~0~~~~2$ & 
3 & 2,2,2\\

$E_6(a_1)$ & $2~~~~0~~~~\overset{\text{\normalsize 0}}{0}~~~~0~~~~2~~~~0~~~~2$ & 
4 & 2,4,6,10\\

$D_5+A_2$ & $0~~~~0~~~~\overset{\text{\normalsize 0}}{0}~~~~2~~~~0~~~~0~~~~2$ & 
3 & 0,2,2\\

$D_6$ & $2~~~~1~~~~\overset{\text{\normalsize 1}}{0}~~~~0~~~~0~~~~1~~~~2$ & 
2 & 2,6\\

$E_6$ & $2~~~~0~~~~\overset{\text{\normalsize 0}}{0}~~~~0~~~~2~~~~2~~~~2$ & 
4 & 2,10,14,22 \\

$D_7(a_2)$ & $1~~~~0~~~~\overset{\text{\normalsize 0}}{1}~~~~0~~~~1~~~~0~~~~1$ & 
3 & 0,2,6\\

$A_7$ & $1~~~~0~~~~\overset{\text{\normalsize 0}}{1}~~~~0~~~~1~~~~1~~~~0$ & 
1 & 2\\

$E_6(a_1)+A_1$ & $1~~~~0~~~~\overset{\text{\normalsize 0}}{1}~~~~0~~~~1~~~~0~~~~2$ & 
3 & 0,2,4\\

$E_7(a_3)$ & $2~~~~0~~~~\overset{\text{\normalsize 0}}{1}~~~~0~~~~1~~~~0~~~~2$ & 
4 & 2,4,6,8\\

$E_8(b_6)$ & $0~~~~0~~~~\overset{\text{\normalsize 0}}{2}~~~~0~~~~0~~~~0~~~~2$ & 
5 & 2,2,2,2,4\\

$D_7(a_1)$ & $2~~~~0~~~~\overset{\text{\normalsize 0}}{0}~~~~2~~~~0~~~~0~~~~2$ & 
4 & 0,2,2,6\\

$E_6+A_1$ & $1~~~~0~~~~\overset{\text{\normalsize 0}}{1}~~~~0~~~~1~~~~2~~~~2$ & 
2 & 2,10\\

$E_7(a_2)$ & $0~~~~1~~~~\overset{\text{\normalsize 1}}{0}~~~~1~~~~0~~~~2~~~~2$ & 
4 & 2,2,6,10\\

$E_8(a_6)$ & $0~~~~0~~~~\overset{\text{\normalsize 0}}{2}~~~~0~~~~0~~~~2~~~~0$ & 
6 & 2,2,2,6,6,6\\

$D_7$ & $2~~~~1~~~~\overset{\text{\normalsize 1}}{0}~~~~1~~~~1~~~~0~~~~1$ & 
2 & 2,10\\

$E_8(b_5)$ & $0~~~~0~~~~\overset{\text{\normalsize 0}}{2}~~~~0~~~~0~~~~2~~~~2$ & 
7 & 2,2,2,2,6,6,10\\

$E_7(a_1)$ & $2~~~~1~~~~\overset{\text{\normalsize 1}}{0}~~~~1~~~~0~~~~2~~~~2$ & 
5 & 2,6,10,14,18\\

$E_8(a_5)$ & $2~~~~0~~~~\overset{\text{\normalsize 0}}{2}~~~~0~~~~0~~~~2~~~~0$ & 
5 & 2,2,2,10,10\\

$E_8(b_4)$ & $2~~~~0~~~~\overset{\text{\normalsize 0}}{2}~~~~0~~~~0~~~~2~~~~2$ & 
5 & 2,2,4,6,10\\

$E_7$ & $2~~~~1~~~~\overset{\text{\normalsize 1}}{0}~~~~1~~~~2~~~~2~~~~2$ & 
4 & 2,10,14,22\\

$E_8(a_4)$ & $2~~~~0~~~~\overset{\text{\normalsize 0}}{2}~~~~0~~~~2~~~~0~~~~2$ & 
6 & 2,4,6,8,10,14\\

$E_8(a_3)$ & $2~~~~0~~~~\overset{\text{\normalsize 0}}{2}~~~~0~~~~2~~~~2~~~~2$ & 
7 & 2,2,8,10,14,16,22\\

$E_8(a_2)$ & $2~~~~2~~~~\overset{\text{\normalsize 2}}{0}~~~~2~~~~0~~~~2~~~~2$ & 
6 & 2,6,10,14,18,22\\

$E_8(a_1)$ & $2~~~~2~~~~\overset{\text{\normalsize 2}}{0}~~~~2~~~~2~~~~2~~~~2$ & 
7 & 2,10,14,18,22,26,34\\

$E_8$ & $2~~~~2~~~~\overset{\text{\normalsize 2}}{2}~~~~2~~~~2~~~~2~~~~2$ & 
8 & 2,14,22,26,34,38,46,58\\

\hline
\end{longtable}

\begin{longtable}{|l|c|c|l|}
\caption{Nilpotent orbits in the Lie algebra of type $F_4$.}\label{tab:F4}
\endfirsthead
\hline
\multicolumn{4}{|l|}{\small\slshape Nilpotent orbits in type $F_4$} \\ 
\hline
\endhead
\hline
\endfoot
\endlastfoot

\hline
label & weighted Dynkin diagram & $\dim \mf{c}_e$ & $h$-weights \\

&  

\begin{picture}(20,7)
  \put(-20,0){\circle{6}}
  \put(0,0){\circle{6}}
  \put(20,0){\circle{6}}
  \put(40,0){\circle{6}}
  \put(-17,0){\line(1,0){14}}
  \put(2,2){\line(1,0){16}}
  \put(2,-2){\line(1,0){16}}
\put(5,-3){$>$}
  \put(23,0){\line(1,0){14}}
\end{picture}

& & \\
\hline

$A_1$ & 1~~~0~~~0~~~0 & 0 &\\

$\widetilde{A}_1$ & 0~~~0~~~0~~~1 &
0 &\\

$A_1+\widetilde{A}_1$ & 0~~~1~~~0~~~0 & 
0 &\\

$A_2$ & 2~~~0~~~0~~~0 & 
1 & 2\\

$\widetilde{A}_2$ & 0~~~0~~~0~~~2 &
1 & 2\\

$A_2+\widetilde{A}_1$ & 0~~~0~~~1~~~0 &
0 & \\

$B_2$ & 2~~~0~~~0~~~1 & 
1 & 2\\

$\widetilde{A}_2+A_1$ & 0~~~1~~~0~~~1 &
1 & 2\\

$C_3(a_1)$ & 1~~~0~~~1~~~0 & 
3 & 2,2,2\\

$F_4(a_3)$ & 0~~~2~~~0~~~0 &
6 & 2,2,2,2,2,2\\

$B_3$ & 2~~~2~~~0~~~0 &
2 & 2,10\\

$C_3$ & 1~~~0~~~1~~~2 &
2 & 2,10\\

$F_4(a_2)$ & 0~~~2~~~0~~~2 &
3 & 2,2,2\\

$F_4(a_1)$ & 2~~~2~~~0~~~2 &
4 & 2,4,6,10\\

$F_4$ & 2~~~2~~~2~~~2 & 
4 & 2,10,14,22\\

\hline
\end{longtable}

\newpage

\begin{longtable}{|l|c|c|l|}
\caption{Nilpotent orbits in the Lie algebra of type $G_2$.}\label{tab:G2}
\endfirsthead
\hline
\multicolumn{4}{|l|}{\small\slshape Nilpotent orbits in type $G_2$} \\
\hline 
\endhead
\hline 
\endfoot
\endlastfoot

\hline
label & weighted Dynkin diagram & $\dim \mf{c}_e$ & $h$-weights \\
&
  \begin{picture}(10,7)
  \put(-10,0){\circle{6}}
  \put(22,0){\circle{6}}
\put(3,-3){$>$}
  \put(-8,-2){\line(1,0){28}}
  \put(-7,0){\line(1,0){26}}
  \put(-8,2){\line(1,0){28}}
\end{picture}

                    & &\\
\hline

$A_1$ & 1~~~~0 & 
1 & 2 \\

$\widetilde{A}_1$ & 0~~~~1 &
0 & \\

$G_2(a_1)$ & 2~~~~0 & 
3 & 2,2,2\\

$G_2$ & 2~~~~2 & 
2 & 2,10\\

\hline

\end{longtable}

\end{document}